\documentclass[a4paper,11pt,english]{article}

\usepackage[french]{babel}

\usepackage[latin1]{inputenc}

\usepackage[T1]{fontenc}

\usepackage{color}

\usepackage[scale={0.7,0.7}]{geometry}

\usepackage{latexsym}
\usepackage[sumlimits,intlimits]{amsmath}
\usepackage{amsthm,amscd,amssymb,amsfonts}
\usepackage{mathrsfs}
\usepackage{url}

\newtheorem{theorem}{Theorem}
\newtheorem{lemma}{Lemma}
\newtheorem{definition}{Definition}
\newtheorem{proposition}{Proposition}
\newtheorem{example}{Example}
\numberwithin{equation}{section}

\usepackage{fancyhdr}
\begin{document}

\title{\textbf{\huge{Regularity Problem for Extremal Vectors}}}
\author{\large Jérôme Verliat$^*$}
\date{}
\maketitle


\thispagestyle{fancy}
\renewcommand{\headrulewidth}{0pt}
\lfoot{\footnotesize{\begin{tabbing} $^*$\= Jérôme Verliat : \= Institut Camille Jordan, Bat. Jean Braconnier, Université Claude Bernard  Lyon 1,\\ \>\> 21 avenue Claude Bernard, 69622 Villeurbanne Cedex, France.\\ \>E-mail address : \url{verliat@math.univ-lyon1.fr} \end{tabbing}}}
\cfoot{}

\vspace{0.8cm}
\begin{center} \textbf{Abstract} \end{center}

\begin{quote}
In this paper, we will use results developed by Ansari and Enflo in the theory of bounded linear operators with dense range. We define two maps, with regards to some parameters, that control surjectivity default of a given operator, and prove analycity for the first one and global continuity for the other one. Minimisation results are also obtained in relation to this study.
\end{quote}


\vspace{0.2cm}
\begin{center}\section*{Introduction}\end{center}

Let $\mathcal{H}$ be a separable infinite dimensional complex Hilbert space and denote by $\mathscr{B}(\mathcal{H})$, the algebra of bounded linear operators on $\mathcal{H}$. Given any operator $T$ of $\mathscr{B}(\mathcal{H})$ with dense range, it is quite natural to estimate it's surjectivity default.

Such a study was firstly proposed by Ansari and Enflo (\cite{ae}) to produce invariant subspaces results by considering the best approximate solutions of the equations $T^ny=x_0$ with $n\in\mathbb{N}$, for the case when $x_0$ is not in the range of $T$. Then, this subject was developed by many authors, for example the contribution of Foias, Jung, Ko and Pearcy (\cite{jkp}, \cite{fjkp2}, \cite{fjkp3}) that made an important breakthrough in the search of hyperinvariant subspaces for some class of operators.

For the sake of continuing this progress, we will work again on the bases of this theory, that is the extremal vectors. Precisely, this paper gives some results concerning minimisation problem in the research of extremal vectors and study their dynamic regarding several parameters.


\vspace{0.2cm}
\begin{center}\item\section{Notations}\end{center}
In the sequel of this paper, $\mathcal{B}(x,r]$ will stand for the closed ball centered on $x$ with radius $r$ and $\partial\mathcal{B}(x,r)$ for its boundary. As $\mathcal{H}$ is a complex Hilbert space, denote by $\langle \cdot | \cdot \rangle$ its complex scalar product, and $[ \cdot | \cdot ]$ the one of the natural real prehilbertian structure of $\mathcal{H}$. For any fixed vector $u$ of $\mathcal{H}$, consider the hyperspace orthogonal to $\textrm{Span}_{\mathbb{R}} (u)$. Then denote by $\mathcal{D}_{u}$ the open half-space delimited by this vector space such that it does not contain $u$ and by $\mathcal{E}_{u}$ its complementary in $\mathcal{H}$.

Let $T$ be an operator of $\mathscr{B}(\mathcal{H})$ with dense range. For such an operator, the vector 0 has a well-known pre-image. Let $x_0$ be a nonzero vector of $\mathcal{H}$ and $\varepsilon$ be any real number such that $0<\varepsilon<\|x_0\|$. Since $T$ has dense range, we can find some image $Ty$ in the ball $\mathcal{B}(x_0,\varepsilon]$. The idea is to pick up one particular vector $y$ such that $Ty$ is in the previous ball. Precisely, Ansari and Enflo consider the $y$ of smallest norm such that $\|Ty-x_0\| \leqslant \varepsilon$ as it is shown in Lemma \ref{lem1}.

\vspace{0.3cm}
\begin{lemma} \label{lem1}
There exists a unique vector $y_0$ that belongs to $\mathcal{B}(x_0,\varepsilon]$ and that satisfies the condition
$$\|y_0\|= \inf \left\{ \|y\|;\|Ty-x_0\| \leqslant \varepsilon \right\}.$$
\end{lemma}

\vspace{0.1cm}
\begin{proof}[Proof]
Set $\mathcal{F}:=\left\{ y\in\mathcal{H}; \| Ty-x_0 \| \leqslant \varepsilon \right\}$. This is a closed convex nonempty subset of $\mathcal{H}$. Indeed, since $T$ has dense range, $\mathcal{F}$ is nonempty. Moreover, we notice that $$T^{-1}\left(\mathcal{B}(x_0,\varepsilon]\right)=\{y\in\mathcal{H};\|Ty-x_0\|\leqslant \varepsilon\}=\mathcal{F}$$
is also a closed set because $T$ is a continuous map. To conclude, if $y_1$ and $y_2$ are two vectors of $\mathcal{H}$ and $t\in[0,1]$, then $\| T(ty_1+(1-t)y_2)-x_0 \| = \| t(Ty_1-x_0)+(1-t)(Ty_2-x_0) \| \leqslant t\varepsilon + (1-t)\varepsilon \leqslant \varepsilon$, so $\mathcal{F}$ is a convex set. Consequently, the projection theorem in Hilbert spaces claims the existence and unicity of $y_0$.
\end{proof}

\vspace{0.1cm}
\begin{definition}
Such a vector is called extremal vector associated with $T$, $x_0$ and $\varepsilon$. It will be denoted by $y_{x_0,\varepsilon}$ in the sequel.
\end{definition}

\vspace{0.1cm}
We have chosen this last notation to stress the dependance of extremal vectors on the two parameters. Indeed, two maps can be defined : either $x_0$ is a fixed nonzero vector, and we define the map
$$\begin{array}{cccl} y_{x_0,\cdot} \ : & ]0;\|x_0\|[ & \longrightarrow & \mathcal{H} \\ & \varepsilon & \longmapsto & y_{x_0,\varepsilon} \end{array} ;$$
or $\varepsilon$ is a fixed nonzero real number, and we define the map
$$\begin{array}{cccl} y_{\cdot,\varepsilon} \ : & \mathcal{H}\backslash\mathcal{B}(0,\varepsilon] & \longrightarrow & \mathcal{H} \\ & x_0 & \longmapsto & y_{x_0,\varepsilon} \end{array} .$$


\vspace{0.2cm}
\begin{center}\section{Minimisation Place Precisions}\end{center}

As mentioned before, to find the extremal vector associated with some $T$, $x_0$ and $\varepsilon$, is a minimisation problem. The computation of $y_{x_0,\varepsilon}$ results from a norm minimisation of a vector $y$ such that $Ty$ describes the whole ball $\mathcal{B}(x_0,\varepsilon]$. In fact, only the minimisation on its boundary is sufficient.

\vspace{0.3cm}
\begin{lemma}[Ansari-Enflo] \label{ae}
The extremal vector associated with $T$, $x_0$ and $\varepsilon$ is the vector $y$ with the smallest norm such that
$$\|Ty-x_0\|=\varepsilon.$$
\end{lemma}

\begin{proof}[Proof]
Let $y_0$ be the extremal vector associated with $T$, $x_0$, $\varepsilon$, and suppose that $\| Ty_0-x_0 \| < \varepsilon$. On the one hand, we can write $\alpha= \varepsilon-\| Ty_0-x_0 \|>0$. On the other hand, there exists a real number $\delta>0$ such that $\mathcal{B}(y_0,\delta]\subset\mathcal{F}$. Indeed, every vector $z$ belonging to $\mathcal{B}(y_0,\delta]$ satisfies : $\|Tz-x_0\| \leqslant \|T(z-y_0)\| + \|Ty_0-x_0\| \leqslant \|T\|\delta + \varepsilon - \alpha \leqslant \varepsilon$ for an appropriate choice of $\delta>0$ provided that $\delta\leqslant \frac{\alpha}{\|T\|}$ holds. Consequently, $(1-\frac{\delta}{2})y_0$ has a smaller norm than $y_0$, which is impossible. The proof is complete.
\end{proof}

We now give a restriction for the minimisation. In that way, we need the following result asserting the collinearity between the two vectors $T^{\star}(Ty_{x_0,\varepsilon}-x_0)$ and $y_0$.

\begin{lemma}[Ansari-Enflo] \label{relation}
There exists a negative number $r_0$ such that
\begin{equation} \label{equ} T^{\star}(Ty_0-x_0)=r_0y_0. \end{equation}
\end{lemma}

\begin{proof}[Proof]
This lemma is based on the following result coming from real prehilbertian structures on a complex Hilbert space $\mathcal{H}$ :
\begin{lemma}
Let $\rho$ be a non negative number. Let $u,v$ be two nonzero vectors of $\mathcal{H}$ such that for all $z\in\mathcal{H}$, $[ u| z] < -\rho\|z\| \Longrightarrow [ v| z] \geqslant \rho\|z\|$ holds.
Then, there exists a negative number $r$ with $|r|>\frac{\rho}{\|u\|}$ such that
$$v=ru.$$
\end{lemma}

\begin{proof}[Proof]
Let $z$ be a vector such that $[ u| z] < -\rho\|z\| \leqslant 0$ : then $z$ belongs to $\mathcal{D}_{u}$. Since $v$ satisfies $[ v| z] \geqslant \rho\|z\| \geqslant 0$, it belongs to $\mathcal{E}_{z}$. Consequently, $v$ belongs to the intersection of all $\mathcal{E}_{z}$ when $z$ describes $\mathcal{D}_{u}$. So there exists a negative number $r$ such that $v=ru$. Then using $z=v$ we get $\|v\| > \rho$, as $\|v\|=|r|\|u\|$, we obtain $|r| > \frac{\rho}{\|u\|}$.
\end{proof}

Therefore, we only need to check that $u=T^{\star}(Ty_0-x_0)$ and $v=y_0$ satisfy this lemma's hypothesis in the particular case when $\rho=0$. Observe that
$$\psi(t)=\|(Ty_0-x_0)+tTz_0\|^2=\|Ty_0-x_0\|^2+2t [ Ty_0-x_0|Tz_0]+t^2\|Tz_0\|^2\geqslant 0$$
and $\psi'(0)=2[ Ty_0-x_0|Tz_0]<0$, thus $\psi$ decreases on $[0,t_0]$ for some $t_0>0$, that is:
$\|T(y_0+tz_0)-x_0\|^2=\psi(t)\leqslant\psi(0)=\|Ty_0-x_0\|^2=\varepsilon^2$.
The minimality of $y_0$ asserts that
$$\|y_0\|^2\leqslant\|y_0+tz_0\|^2=\|y_0\|^2+2t[ y_0|z_0]+t^2\|z_0\|^2,$$
consequently $\varphi(t):=\|y_0+tz_0\|^2$ defines a decreasing map on $[0,t_0]$, implying $\varphi'(0)\geqslant0$, \textit{i.e.}, $[y_0|z_0]\geqslant0$.
\end{proof}

This result, proved in \cite{ae}, provides the following new reduction of the minimisation place. We are only satisfied with the minimisation on a kind of cap of the sphere $\partial \mathcal{B}(x_0,\varepsilon)$. More precisely:

\vspace{0.3cm}
\begin{proposition} \label{proplieu}
The extremal vector $y_{x_0,\varepsilon}$ is the vector of the smallest norm such that $Ty_{x_0,\varepsilon}$ belongs to
$$\mathcal{V}_{x_0,\varepsilon}:=\partial\mathcal{B}(x_0,\varepsilon)\cap\mathcal{B}(0,\sqrt{\|x_0\|^2-\varepsilon^2}].$$
\end{proposition}

\begin{proof}[Proof]
Lemma \ref{relation} implies the existence of a negative number $r_0$ such that $T^{\star}(Ty_0-x_0)=r_0y_0$. This entails $[ Ty_0-x_0|Ty_0] = [ T^{\star}(Ty_0-x_0)|y_0]=r_0\|y_0\|^2 <0$, meaning that the angle formed by vectors $Ty_0-x_0$ and $Ty_0$ is obtuse. For a better understanding, we can put forward the following picture, explaining the situation in the two dimension case:
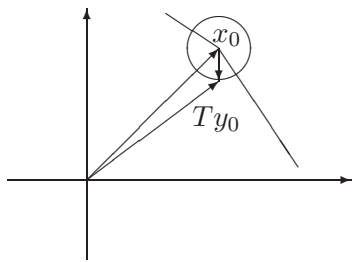
\begin{center}
    \begin{picture}(80,100)(-50,50)
          \put(-40,80){\vector(1,0){100}}
          \put(-40,80){\vector(0,1){65}}
          \put(-40,80){\vector(1,1){50}}
          \put(-40,80){\vector(4,3){50}}
          \put(10,130){\vector(0,-1){12.5}}
          \put(-40,80){\line(-1,0){30}}
          \put(-40,80){\line(0,-1){30}}
          \put(13,134){\makebox(0,0){$x_0$}}
          \put(9,103){\makebox(0,0){$Ty_0$}}
          \put(10,130){\circle{25.5}}
          \put(10,130){\line(-3,2){20}}
          \put(10,130){\line(2,-3){30}}
    \end{picture}
\end{center}
Vectors $y$, such that $Ty$ are situated on the right hand side of the two half straight lines drawn from $x_0$ endpoint, are not suitable. Indeed, in this case $Ty-x_0$ and $Ty$ form an obtuse angle. Thus, such a vector cannot be the extremal vector $y_0$ and it will not be considered in the minimisation. The infinite dimension case is similar: the only vectors $y$ that play a rôle in the minimisation are the ones such that $Ty$ are situated on one side of a half cone given by the straight lines through $x_0$ and orthogonal to tangents through 0 of the circle centered on $x_0$ with radius $\varepsilon$. The two dimension case can be helpful :
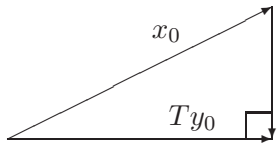
\begin{center}
    \begin{picture}(100,60)(-50,55)
          \put(-50,60){\vector(2,1){100}}
          \put(-50,60){\vector(1,0){100}}
          \put(50,110){\vector(0,-1){50}}
          \put(50,70){\line(-1,0){10}}
          \put(40,60){\line(0,1){10}}
          \put(10,100){\makebox(0,0){$x_0$}}
          \put(20,68){\makebox(0,0){$Ty_0$}}
    \end{picture}
\end{center}
Use Pythagorean theorem to observe that those vectors norms are $\sqrt{\|x_0\|^2-\varepsilon^2}$.
\end{proof}


\vspace{0.2cm}
\begin{center} \section{Regularity with regard to $\varepsilon$} \end{center}

Fix now a vector $x_0\in\mathcal{H}$ and consider extremal vectors associated with $T$, $x_0$ and $\varepsilon$ where $\varepsilon\in]0,\|x_0\|[$. We focus now on the map defined in Section 1 :
$$\begin{array}{cccl} y_{x_0,\cdot} \ : & ]0;\|x_0\|[ & \longrightarrow & \mathcal{H} \\ & \varepsilon & \longmapsto & y_{x_0,\varepsilon} \end{array}.$$
The next result, proved in \cite{ae}, holds :

\vspace{0.3cm}
\begin{proposition}[Ansari-Enflo] \label{propenflo}
The map $\varepsilon \mapsto y_0$ is analytic over $]0,\|x_0\|[$.
\end{proposition}

\begin{proof}[Proof]
Observe that it suffices to prove the analycity of $\varepsilon\mapsto r_{\varepsilon}$, where $r_{\varepsilon}$ is the negative number given by Lemma \ref{relation}. Indeed, Equation \ref{equ} implies $(r_{\varepsilon}I-T^{\star}T)y_0=-T^{\star}x_0$, and since $r_{\varepsilon}I-T^{\star}T$ is invertible (because for all $r<0$, $T^{\star}T-rI\geqslant -rI>0$), $y_0=-(r_{\varepsilon}I-T^{\star}T)^{-1}T^{\star}x_0$ holds. Focus now on the converse map of $\varepsilon\mapsto r_{\varepsilon}$, that is $r_{\varepsilon}\mapsto \varepsilon$, which is well defined since $y_0$ is unique. As
$\varepsilon^2=\|Ty_0-x_0\|^2=\|-T(r_{\varepsilon}I-T^{\star}T)^{-1}T^{\star}x_0-x_0\|^2$,
this entails the analycity of the map $r_{\varepsilon}\mapsto \varepsilon^2$, which is moreover injective. Consequently, $r_{\varepsilon}\mapsto \varepsilon$ is analytic.
\end{proof}


\vspace{0.2cm}
\begin{center} \section{Regularity with regard to $x_0$} \end{center}

In this section, we fix the real number $\varepsilon$ to be positive, and study extremal vectors associated with $T$, $\varepsilon$ and $x_0$ where $x_0$ is a vector in $\mathcal{H}\backslash \mathcal{B}(0,\varepsilon]$. In a first part, we will obtain the continuity of the map defined in Section 1 :
$$\begin{array}{cccl} y_{\cdot,\varepsilon} \ : & \mathcal{H}\backslash\mathcal{B}(0,\varepsilon] & \longrightarrow & \mathcal{H} \\ & x_0 & \longmapsto & y_{x_0,\varepsilon} \end{array} .$$
To prove this result, we firstly consider the map $x_0\mapsto \|y_{T,x_0,\varepsilon}\|$ :

\vspace{0.3cm}
\begin{proposition} \label{contennorm}
The map $x_0 \mapsto \|y_{x_0,\varepsilon}\|$ is continuous on the open set $\mathcal{H}\backslash\overline{\mathcal{B}(x_0,\varepsilon)}$.
\end{proposition}

\vspace{0.1cm}
\begin{proof}[Proof]
Let $\mu$ be a positive number. The analycity of $\varepsilon\mapsto y_{x_0,\varepsilon}$ implies the continuity of $\varepsilon\mapsto \|y_{x_0,\varepsilon}\|$. Thus, there exists $\nu>0$ such that $\displaystyle \left| \|y_{x_0,\varepsilon \pm \nu}\| - \|y_{x_0,\varepsilon}\| \right| < \mu$. Let $x\in\mathcal{B}(x_0,\nu]$. Then $\mathcal{B}(x_0,\varepsilon-\nu] \subset \mathcal{B}(x,\varepsilon] \subset \mathcal{B}(x_0,\varepsilon+\nu]$. Consequently, minimality's definition involves $\|y_{x_0,\varepsilon + \nu}\| \leqslant \|y_{x,\varepsilon}\| \leqslant \|y_{x_0,\varepsilon - \nu}\|$. So, vectors $x$ of the ball $\mathcal{B}(x_0,\nu]$ satisfy
$$\|y_{x_0,\varepsilon}\|-\mu \leqslant \|y_{x,\varepsilon}\| \leqslant \|y_{x_0,\varepsilon}\|+\mu$$
which completes the proof.
\end{proof}

This result will help us to prove the global continuity of the map $x_0 \mapsto y_{x_0,\varepsilon}$. The proof use Proposition \ref{contennorm}, the Pythagorean theorem and a Banach space theorem, concerning a decreasing intersection of closed sets whose diameters converge to 0.

\vspace{0.3cm}
\begin{theorem}
With the same hypotheses as in Proposition \ref{contennorm}, the map $x_0\mapsto y_0$ is continuous at all point of the open set $\mathcal{H}\backslash\mathcal{B}(0,\varepsilon]$.
\end{theorem}

\vspace{0.1cm}
\begin{proof}[Proof]
We restrict the study of the continuity at the vector $x_0$. Let $x_0^{(n)}$ be a sequence converging in $\mathcal{H}$ to $x_0$. We set $y_0^{(n)}:=y_{x_0^{(n)},\varepsilon}$. Notice that for all positive number $\nu$, there exists an integer $N$ such that
$$\forall n\geqslant N, \quad y_0^{(n)}\in \left\{ y\in\mathcal{H} ; \|Ty-x_0\| \leqslant \varepsilon+\nu \right\}.$$
Furthermore, the vector $y_{x_0,\varepsilon}$ belongs to this set. For this reason a decreasing sequence $(\nu_n)_n$ of positive numbers converging to 0, and such that for all integer $k\geqslant n$, $y_0^{(k)} \in \left\{ y\in\mathcal{H} ; \|Ty-x_0\| \leqslant \varepsilon+\nu_n \right\}$ can be constructed. Set :
$$F_n:=\left\{ y\in\mathcal{H} ; \|Ty-x_0\| \leqslant \varepsilon+\nu_n \textrm{ and } \|y\| \leqslant \sup_{k\geqslant n}\left\{\|y_0^{(k)}\|,\|y_{x_0,\varepsilon}\|\right\} \right\}.$$
This is a closed convex set, and $y_0^{(n)}\in F_n$. With the Pythagorean theorem, its diameter is bounded by the real number
$$2\times\sqrt{\sup_{k\geqslant n} \left\{\|y_0^{(k)}\|,\|y_{x_0,\varepsilon}\|\right\}^2-\|y_{x_0,\varepsilon+\nu_n}\|^2}\ .$$
Therefore $(F_n)_{n\in\mathbb{N}}$ is a decreasing sequence of closed sets whose diameters converge to 0. Completeness of the Hilbert space $\mathcal{H}$ involves
$$\bigcap_{n\in\mathbb{N}} F_n = \left\{ y_{x_0,\varepsilon} \right\}.$$
Consequently the sequence $(y_0^{(n)})$ converges to $y_{x_0,\varepsilon}$.
\end{proof}

In this part, we study a particular case of the previous application. Here, $T$ and $\varepsilon$ are fixed as before, and $x_0$ is a fixed vector in $\mathcal{H}\backslash\mathcal{B}(0,\varepsilon]$. Study now the application
$$\left\{ \begin{array}{ccl} \displaystyle\left]\frac{\varepsilon}{\|x_0\|},+\infty\right[ & \longrightarrow & \mathbb{R}\\ t & \longmapsto & \|y_{tx_0,\varepsilon}\| \end{array} \right. .$$

\begin{proposition} \label{propcroissance}
The map $t\longmapsto \|y_{tx_0,\varepsilon}\|$ is increasing on $\displaystyle\left]\frac{\varepsilon}{\|x_0\|},+\infty\right[$.
\end{proposition}

\begin{proof}[Proof]
For more simplicity, consider the map $t\mapsto \|y_{(1+t)x_0,\varepsilon}\|$ instead of the one of that Proposition.

Firstly, observe that
$$\mathcal{V}_{(1+t)x_0,\varepsilon}\subset\mathcal{B}(x_0,\varepsilon]$$
for all $t\in\mathbb{R}_+^{\star}$ sufficiently small, for example since $t$ belongs to $\left[0,\frac{\varepsilon^2}{\|x_0\|^2}\right]$. Indeed, if $x'\in\mathcal{V}_{(1+t)x_0,\varepsilon}$, then we obtain $\|(1+t)x_0-x'\|=\varepsilon$ and $\|x'\|^2\leqslant (1+t)^2\|x_0\|^2-\varepsilon^2$. Writing $x'=(1+t)x_0+k$ for $k\in\mathcal{H}$, we can develop : $\|x'\|^2=(1+t)^2\|x_0\|^2+2(1+t)[x_0|k]+\varepsilon^2\leqslant(1+t)^2\|x_0\|^2-\varepsilon^2$, so
$$\varepsilon^2+(1+t)[x_0|k]\leqslant0$$
holds. Furthermore $x'-x_0=tx_0+k$, so $\|x'-x_0\|\leqslant\varepsilon$ is equivalent to $t\|x_0\|^2+2[x_0|k]\leqslant0$.
Since $t\|x_0\|^2+2[x_0|k]\leqslant t\|x_0\|^2-2\frac{\varepsilon^2}{1+t}$, it is sufficient to find non negative numbers $t$ such that
$$(1+t)t\|x_0\|^2- 2\varepsilon^2\leqslant0.$$
A short study of this polynomial shows it vanishes twice, one of its zero is negative and the other one is positive, equals to $\frac{1}{2}\left(\frac{\sqrt{\|x_0\|^2+8\varepsilon^2}}{\|x_0\|}-1\right)$. We only have to check that the previous bound is less or equal to this one:
$$\frac{\sqrt{\|x_0\|^2+8\varepsilon^2}}{2\|x_0\|}-\frac{1}{2}=\frac{4\varepsilon^2}{\|x_0\|(\sqrt{\|x_0\|^2+8\varepsilon^2}+\|x_0\|)}\geqslant\frac{\varepsilon^2}{\|x_0\|^2}$$
because $\varepsilon<\|x_0\|$.

Consequently, using minimality argument of Proposition \ref{proplieu}, we obtain $\|y_{x_0,\varepsilon}\| \leqslant \|y_{(1+t)x_0,\varepsilon}\|$. Since $x_0$ is arbitrary, $\|y_{(1+t_1)x_0,\varepsilon}\| \leqslant \|y_{(1+t_2)x_0,\varepsilon}\|$ for all $t_1\leqslant t_2$, and the map is increasing on $\left[0,\frac{\varepsilon^2}{\|x_0\|^2}\right]$. The same idea applied to $x_1=(1+\frac{\varepsilon^2}{\|x_0\|^2})x_0$ instead of $x_0$ shows that $t\mapsto \|y_{(1+t)x_1,\varepsilon}\|$ is increasing on $\left[0,\frac{\varepsilon^2}{\|x_1\|^2}\right]$. Applying this argument again, every map $t\mapsto \|y_{(1+t)x_n,\varepsilon}\|$ is increasing on $\left[0,\frac{\varepsilon^2}{\|x_n\|^2}\right]$, where the sequence $(x_n)_n$ is defined by $x_0$ and $x_{n+1}=(1+\frac{\varepsilon^2}{\|x_n\|^2})x_n$. We get $\|x_{n+1}\|=(1+\frac{\varepsilon^2}{\|x_n\|^2})\|x_n\|=\|x_n\|+\frac{\varepsilon^2}{\|x_n\|}$. As the sequence $(\|x_n\|)_n$ is increasing, and the map $f:x\mapsto x+\frac{\varepsilon^2}{x}$ has no fixed point, $(\|x_n\|)_n$ tends to infinity as $n$ tends to infinity. Therefore the map $t\mapsto \|y_{(1+t)x_0,\varepsilon}\|$ is increasing on $[0,+\infty[$. Taking another $x_0$ in the same direction but with a norm closed to $\varepsilon$, we obtain that this map is increasing on $\left]\frac{\varepsilon}{\|x_0\|}-1,+\infty\right[$, and the proof is complete.
\end{proof}

\paragraph{Counter-example in the two dimensional case.}
We show now that Proposition~\ref{propcroissance} cannot be generalised to all directions thanks to a counter-example in the particular case of dimension two. Let $\mathcal{H}$ be the two dimensional complex Hilbert space $\mathbb{C}^2$.

\begin{example}
Consider $T=I_{\mathbb{C}^2}$ be the identity operator of $\mathscr{B}(\mathbb{C}^2)$. Let $x_0=\begin{pmatrix} 2 \\ -2 \end{pmatrix}$ and $u=\begin{pmatrix} 0 \\ 2 \end{pmatrix}$ be two vectors of $\mathbb{C}^2$. Let $\varepsilon=1$. Then, the map $t\mapsto \|y_{x_0+tu,\varepsilon}\|$ is decreasing on $[0;1]$ and increasing on $[1;+\infty[$.
\end{example}

\begin{proof}[Proof]
We know that the extremal vector $y_{x_0+tu,\varepsilon}$ belongs to $T^{-1}(\partial\mathcal{B}(x_0+tu,\varepsilon))$, which is determined by the equation $\|Ty-x_0-tu\|_2=\varepsilon$. Writing $y=\begin{pmatrix} \lambda \\ \mu \end{pmatrix}$, it means that $\lambda$ and $\mu$ satisfy $(\lambda-2)^2+(\mu+2-2t)^2=1$. This is the  equation of a circle $\mathscr{C}$ centered on $C=\begin{pmatrix} 2 \\ 2t-2 \end{pmatrix}$ with radius $r=1$. Therefore $$\|y_{x_0+tu,\varepsilon}\|=\textrm{d}(O,\mathscr{C})=OC-r=\sqrt{4+4(t-1)^2}-1=2\sqrt{(t-1)^2+1}-1$$
where d stands for the distance in $\mathbb{C}^2$. A short study of the function $t\mapsto 2\sqrt{(t-1)^2+1}-1$ completes the proof.
\end{proof}


\vspace{0.5cm}
\underline{Acknowledgements}: The author wishes to thank University of Lyon, France, where part of this work was presented; G. Cassier for its help and its good advices, and P. Leroux as well for advices in the translation of this work.


\vspace{0.2cm}
\renewcommand{\refname}{Bibliography}

\end{document}